\documentclass[12pt]{amsart}%
\usepackage{graphicx}
\usepackage{amscd}
\usepackage{amsmath}
\usepackage{amsfonts}
\usepackage{amssymb}%
\providecommand{\U}[1]{\protect\rule{.1in}{.1in}}
\theoremstyle{plain}

\theoremstyle{definition}

\newtheorem{remark}{Remark}[section]

\numberwithin{equation}{section}
\numberwithin{theorem}{section}
\begin{document}
\title[Orthogonal polynomials]{Note on the $X_{1}$-\textbf{Laguerre} orthogonal polynomials.}
\author{W.N. Everitt}
\address{W.N. Everitt, School of Mathematics, University of Birmingham, Edgbaston,
Birmingham B15 2TT, England, UK}
\email{w.n.everitt@bham.ac.uk}
\date{21 November 2008 (File C:%
$\backslash$%
Swp55%
$\backslash$%
Docs%
$\backslash$%
milson22.tex)}

\begin{abstract}
This note supplements the results in the paper on $X_{1}$-\textbf{Laguerre}
orthogonal polynomials written by David G\'{o}mez-Ullate, Niky Kamran and
Robert Milson.

\end{abstract}
\maketitle

\section{Introduction\label{sec1}}

This note reports on, the $X_{1}$-\textbf{Laguerre polynomials}, one of the
two new sets of orthogonal polynomials considered in the papers \cite{G-UKM}
and \cite{G-UKM1}, written by David G\'{o}mez-Ullate, Niky Kamran and Robert
Milson. The other set is named the $X_{1}$-\textbf{Jacobi polynomials} and is
discussed, in similar terms, in the note \cite{WNE2}.

These two papers are remarkable and invite comments on the results therein
which have yielded new examples of Sturm-Liouville differential equations and
their associated differential operators.

The two sets of these orthogonal polynomials are distinguished by:

\begin{itemize}
\item[$(i)$] Each set of polynomials is of the form $\{P_{n}(x):x\in
\mathbb{R}\ $and$\ n\in\mathbb{N}\equiv\{1,2.3.\ldots\}\}$ with$\ \deg
(P_{n})=n;$ that is there is no polynomial of degree $0.$

\item[$(ii)$] Each set is orthogonal and complete in a weighted Hilbert
function space.

\item[$(iii)$] Each set is generated as a set of eigenvectors from a
self-adjoint Sturm-Liouville differential operator.
\end{itemize}

\section{$X_{1}$-\textbf{Laguerre polynomials\label{sec2}}}

These polynomials and the associated differential equation are detailed in
\cite[Section 2]{G-UKM}.

In \cite[Section 2, (21)]{G-UKM} the second-order linear differential equation
concerned is given as%
\begin{equation}
-xy^{\prime\prime}(x)+\left(  \dfrac{x-k}{x+k}\right)  \left(
(x+k+1)y^{\prime}(x)-y(x)\right)  =\lambda y(x)\ \text{for all}\ x\in
(0,\infty) \label{eq2.0}%
\end{equation}
where the parameter $\lambda\in\mathbb{C}$ plays the role of a spectral
parameter for the differential operators defined below, and the parameter
$k\in(0,\infty).$

This equation (\ref{eq2.0}) is not a Sturm-Liouville differential equation;
such equations take the form, in this case taking the interval to be
$(o,\infty),$%
\begin{equation}
-(p(x)y^{\prime}(x))^{\prime}+q(x)y(x)=\lambda w(x)y(x)\ \text{for all}%
\ x\in(0,\infty), \label{eq2.0a}%
\end{equation}
but can be transformed into this form on using the information in
\cite[Section 2]{G-UKM}. In particular let the coefficients $p_{k},q_{k}%
,w_{k}$ be defined as follows;

\begin{itemize}
\item[$(i)$] $p_{k},q_{k},w_{k}:(0,\infty)\rightarrow\mathbb{R}$

\item[$(ii)$]
\begin{equation}
p_{k}(x):=\frac{x^{k}}{(x+k)^{2}}\exp(-x)\ \text{for all }x\in(0,\infty)
\label{eq2.1}%
\end{equation}

\item[$(iii)$]
\begin{equation}
q_{k}(x):=-\frac{(x-k)x^{k}}{(x+k)^{3}}\exp(-x)~\text{for all }x\in(0,\infty)
\label{eq2.2}%
\end{equation}

\end{itemize}

\begin{equation}
w_{k}(x):=\frac{x^{k}}{(x+k)^{2}}\exp(-x)\ \text{for all }x\in(0,\infty).
\label{eq2.3}%
\end{equation}

Let the Sturm-Liouville differential expression $M_{k}$ have the domain%
\begin{equation}
D(M_{k}):=\{f:(0,\infty)\rightarrow\mathbb{C}:f^{(r)}\in AC_{\text{loc}%
}(0,\infty)\ \text{for}\ r=0,1\} \label{eq2.4}%
\end{equation}
and be defined by, for all $f\in D(M_{k}),$%
\begin{equation}
M_{k}[f](x):=-(p_{k}(x)f^{\prime}(x))^{\prime}+q_{k}(x)f(x)\ \text{for almost
all}\ x\in(0,\infty). \label{eq2.5}%
\end{equation}

Now define the Sturm-Liouville differential equation by, for all
$k\in(0,\infty),$%
\begin{equation}
M_{k}[y](x)=\lambda w_{k}(x)y(x)\ \text{for all}\ x\in(0,\infty) \label{eq2.6}%
\end{equation}
where $\lambda\in\mathbb{C}$ is a complex valued spectral parameter.

For an account of Sturm-Liouville theory of differential operators and
equations, see \cite[Sections 2 to 6]{WNE}.

It is important to notice that the differential equation (\ref{eq2.6}) is
equivalent to, and is derived from the differential equation (\ref{eq2.0}),
see again \cite[Section 2.2, (22a)]{G-UKM}.

The differential equation (\ref{eq2.6}) is to be studied in the Hilbert
function space $L^{2}((0,\infty);w_{k}).$

The symplectic form for $M_{k}$ is defined by, for all $k\in(0,\infty)$ and
for all $f,g\in D(M_{k}),$%
\begin{equation}
\lbrack f,g]_{k}(x):=f(x)(p_{k}\overline{g}^{\prime})(x)-(p_{k}f^{\prime
})(x)\overline{g}(x)\ \text{for all}\ x\in(0,\infty). \label{eq2.7}%
\end{equation}

The maximal operator $T_{k,1}$ is defined by, for all $k\in(0,\infty),$%
\begin{equation}
\left\{
\begin{array}
[c]{ll}%
(i) & T_{k,1}:D(T_{1,k})\subset L^{2}((0,\infty);w_{k})\rightarrow
L^{2}((0,\infty);w_{k})\\
(ii) & D(T_{k,1}):=\{f\in D(M_{k}):f,w^{-1}M_{k}[f]\in L^{2}((0,\infty
);w_{k})\}\\
(iii) & T_{k,1}f:=w^{-1}M_{k}[f]\ \text{for all}\ f\in D(T_{1,k}).
\end{array}
\right.  \label{eq2.8}%
\end{equation}

All self-adjoint differential operators in $L^{2}((0,\infty);w_{k})$ generated
by $M_{k}$ are given by restrictions of the maximal operator $T_{k,1};$ these
restrictions are determined by placing boundary conditions at the endpoints
$0$ and $\infty,$ on the elements of $D(T_{k,1}).$ The number and type of
boundary conditions depends upon the endpoint classification of $M_{k}$ in
$L^{2}((0,\infty);w_{k});$ see \cite[Section 5]{WNE}.

For the endpoint classification of the differential expression $M_{k}$ in
$L^{2}((0,\infty);w_{k})$ we have the results, see again \cite[Section 5]{WNE};

\begin{itemize}
\item[$(i)$] At $0^{+}$ the classification is:%
\begin{equation}%
\begin{tabular}
[c]{|l|l|}\hline
For all $k\in(0,3]$ & limit-circle non-oscillatory\\\hline
& \\\hline
For all $k\in(3,\infty)$ & limit-point.\\\hline
\end{tabular}
\ \ \label{eq2.9}%
\end{equation}

\item[$(ii)$] At $+\infty$ the classification is:%
\begin{equation}%
\begin{tabular}
[c]{|l|l|}\hline
For all $k\in(0,\infty)$ & limit point.\\\hline
\end{tabular}
\ \ \label{eq2.10}%
\end{equation}

\end{itemize}

To establish these properties we have the following results:

\begin{enumerate}
\item For $\lambda=0$ the function%
\begin{equation}
\varphi_{1}(x):=x+k+1\ \text{for all}\ x\in\lbrack0,\infty), \label{eq2.11}%
\end{equation}
is a solution of the differential equation (\ref{eq2.6}), for all
$k\in(0,\infty);$ see \cite[Section 2, (14)]{G-UKM}.

\item We have $\varphi_{1}\in L^{2}((0,\infty);w_{k})$ for all $k\in
(0,\infty).$

\item For $\lambda=0$ the function%
\begin{equation}
\varphi_{2}(x):=\varphi_{1}(x)\int_{1}^{x}\frac{1}{\varphi_{1}^{2}(t)p_{k}%
(t)}~dt~\text{for all}\ x\in(0,\infty), \label{eq2.12}%
\end{equation}
is a solution of the differential equation (\ref{eq2.6}), for all
$k\in(0,\infty);$ $\varphi_{2}$ is independent of $\varphi_{1}.$

\item Asymptotic analysis shows that%
\begin{equation}%
\begin{tabular}
[c]{|l|}\hline
$\varphi_{1}\in L^{2}((0,\infty);w_{k})$ for all $k\in(0,\infty)$\\\hline
\end{tabular}
\ \ \ \ \label{eq2.13}%
\end{equation}
and%
\begin{equation}%
\begin{tabular}
[c]{|l|}\hline
$\varphi_{2}\notin L^{2}([1,\infty);w_{k})$ for all $k\in(0,\infty)$\\\hline
\\\hline
$\varphi_{2}\in L^{2}((0,1];w_{k})$ for all $k\in(0,3]$\\\hline
\\\hline
$\varphi_{2}\notin L^{2}((0,1];w_{k})$ for all $k\in(3,\infty).$\\\hline
\end{tabular}
\ \ \ \ \label{eq2.14}%
\end{equation}

\end{enumerate}

The endpoint classifications (\ref{eq2.9}) and (\ref{eq2.10}) follow from the
results items 1 to 4 above; see \cite[Section 5.]{WNE}.

We can now define the restriction $A_{k}$ of the maximal operator $T_{k,1},$
see (\ref{eq2.8}), which is self-adjoint in the Hilbert function space
$L^{2}((0,\infty);w_{k}),$ and which has the $X_{1}$-Laguerre polynomials as
eigenvectors. To obtain this result it is essential:

\begin{itemize}
\item[$(i)$] To apply the general theory of such restrictions as given in the
Naimark text \cite[Chapter V, Sections 17 and 18]{MAN}.

\item[$(ii)$] To apply the detailed results on the properties of the $X_{1}%
$-Laguerre polynomials given in \cite[Section 2]{G-UKM}.
\end{itemize}

At any limit-point endpoint no boundary condition is required; at the
limit-circle endpoint $0^{+}$ the boundary condition for any $f\in D(T_{k,1})$
takes the form%
\begin{equation}
\lim_{x\rightarrow0+}\left[  f,\alpha_{1}\varphi_{1}+\alpha_{2}\varphi
_{2}\right]  (x)=0, \label{eq2.15}%
\end{equation}
where $\alpha_{1},\alpha_{2}\in\mathbb{R}.$ Since the $X_{1}$-Laguerre
polynomials are to be in the domain of the operator $A_{k}$ we take
$\alpha_{1}=1$ and $\alpha_{2}=0.$

Thus the domain $D(A_{k})$ of our self-adjoint operator $A_{k}$ restriction of
the maximal operator $T_{k}$ is defined as follows:

\begin{itemize}
\item[$(i)$] For $k\in(0,3]$%
\begin{equation}
D(A_{k}):=\{f\in D(T_{k,1}):\lim_{x\rightarrow0+}\left[  f,\varphi_{1}\right]
(x)=0\} \label{eq2.16}%
\end{equation}
and%
\begin{equation}
A_{k}f:=w_{k}^{-1}M_{k}[f]\ \text{for all}\ f\in D(A_{k}). \label{eq2.17}%
\end{equation}

\item[$(ii)$] For $k\in(3,\infty)$%
\begin{equation}
D(A_{k}):=D(T_{k,1}) \label{eq2.18}%
\end{equation}
and%
\begin{equation}
A_{k}f:=w_{k}^{-1}M_{k}[f]\ \text{for all}\ f\in D(A_{k}). \label{eq2.19}%
\end{equation}

\end{itemize}

The spectrum and eigenvectors of $A_{\alpha,\beta}$ can be obtained from the
results given in \cite[Section 2]{G-UKM}. The spectrum of $A_{\alpha,\beta}$
contains the sequence $\{\lambda_{n}=n:n\in\mathbb{N}_{0}\};$ the eigenvectors
are given by $\{\hat{L}_{n+1}^{(k)}:n\in\mathbb{N}_{0}\},$ the $X_{1}%
$-\textbf{Laguerre} orthogonal polynomials.

\begin{remark}
\label{rem2.1}

\begin{itemize}
\item[$(i)$] The notation $\lambda_{n}=n$ for all $n\in\mathbb{N}_{0}$ makes
good comparison with the eigenvalue notation for the classical Laguerre
polynomials; this sequence is independent of the parameter $k\in(0,\infty).$

\item[$(ii)$] We note that $\hat{L}_{n+1}^{(k)}$ is a polynomial of degree
$n+1$ for all $n\in\mathbb{N}_{0}$ and all $k\in(0,\infty).$

\item[$(iii)$] Note that for $k\in(0,3],$ when the limit-circle condition
holds at $0^{+},$ it is essential to check that the polynomials $\{\hat
{L}_{n+1}^{(k)}\}$ all satisfy the boundary condition at $0^{+}$ as required
in $(\ref{eq2.16}),$ \textit{i.e.}%
\begin{equation}
\lim_{x\rightarrow0+}\left[  \hat{L}_{n+1}^{(k)},\varphi_{1}\right]
(x)=0\ \text{for all}\ n\in\mathbb{N}_{0.} \label{eq2.21}%
\end{equation}
This result follows since, using $(\ref{eq2.11}),$%
\begin{align*}
\left[  \hat{L}_{n+1}^{(k)},\varphi_{1}\right]  (x)  &  =p_{k}(x)\left[
\hat{L}_{n+1}^{(k)}(x)\varphi_{1}^{\prime}(x)-\hat{L}_{n+1}^{(k)\prime}%
\varphi_{1}(x)\right] \\
&  =\dfrac{x^{k}}{(x+k)^{2}}\exp(-x)\left[  \hat{L}_{n+1}^{(k)}-\hat{L}%
_{n+1}^{(k)\prime}(x+k+1)\right] \\
&  =\mathcal{O}(x^{k})\ \text{as}\ x\rightarrow0^{+}.
\end{align*}

\end{itemize}
\end{remark}

It is shown in \cite[Section 3, Proposition 3.3]{G-UKM} that the sequence of
polynomials
\[
\left\{  \hat{L}_{n+1}^{(k)}:n\in\mathbb{N}_{0}\right\}
\]
is orthogonal and dense in the space $L^{2}((0,\infty);w_{k}),$ for all
$k\in(0,\infty).$ This result implies that for all $k\in(0,\infty)$ the
spectrum of the operator $A_{k}$ consists entirely of the sequence of
eigenvalues $\{\lambda_{n}:n\in\mathbb{N}_{0}\};$ from the spectral theorem
for self-adjoint operators in Hilbert space it follows that no other point on
the real line $\mathbb{R}$ can belong to the spectrum of $A_{k}.$

\begin{remark}
\label{rem2.2}It is to be noted that whilst the Hilbert space theory as given
in \cite{WNE} and \cite{MAN} provides a precise definition of the self-adjoint
operator $A_{k},$ the information about the particular spectral properties of
$A_{k}$ are to be deduced from the classical analysis results in \cite{G-UKM}.
Without these results it would be very difficult to deduce the spectral
properties of the self-adjoint operator $A_{k},$ as defined above, in the
Hilbert function space $L^{2}((0,\infty);w_{k}).$
\end{remark}

\end{document}